\documentclass{amsart}
\usepackage{amssymb}
\usepackage{enumerate}
\usepackage{url}
\usepackage{verbatim}

\newcommand\N{\mathbb N}
\newcommand\Z{\mathbb Z}

\newcommand\R{\mathbb R}
\newcommand\C{\mathbb C}

\newcommand\ph\varphi
\newcommand\ps\psi
\newcommand\ep\varepsilon
\newcommand\rh\varrho
\newcommand\al\alpha
\newcommand\be\beta
\newcommand\ga\gamma
\newcommand\om\omega
\newcommand\ta\tau
\renewcommand\th\vartheta
\newcommand\de\delta
\newcommand\ze\zeta
\newcommand\ch\chi
\newcommand\et\eta
\newcommand\io\iota
\newcommand\la\lambda
\newcommand\si\sigma

\newcommand\Ga\Gamma
\newcommand\De\Delta
\newcommand\Th\Theta
\newcommand\La\Lambda
\newcommand\Si\Sigma
\newcommand\Ph\Phi
\newcommand\Ps\Psi
\newcommand\Om\Omega

\newtheorem{theorem}{Theorem}
\newtheorem{openproblem}[theorem]{Open Problem}
\newtheorem{lemma}[theorem]{Lemma}
\newtheorem{proposition}[theorem]{Proposition}
\newtheorem{corollary}[theorem]{Corollary}
\theoremstyle{definition}
\newtheorem{definition}[theorem]{Definition}

\theoremstyle{remark}
\newtheorem{remark}[theorem]{Remark}
\newtheorem{example}[theorem]{Example}

\newcommand\x{{\bar X}}
\newcommand\rx{{\R[\x]}}
\newcommand\nn{\R_{\ge 0}}
\newcommand\pos{\R_{>0}}

\newcommand{\ran}{R_\infty}
\newcommand{\ptent}[1]{S(\nabla#1)}
\newcommand{\tent}[2]{S(\nabla#1,#2)}
\newcommand{\sos}{\text{sos}}

\begin{document}
\title[Global optimization of polynomials]
{Global optimization of polynomials using gradient tentacles and
sums of squares}
\author{Markus Schweighofer}
\thanks{Supported by the DFG grant ``Barrieren''.}
\address{Universit\"at Konstanz\\
         Fachbereich Mathematik und Statistik\\
         78457 Konstanz\\
         Allemagne}
\email{Markus.Schweighofer@uni-konstanz.de}
\keywords{global optimization, polynomial, preorder, sum of squares,
semidefinite programming}
\subjclass{Primary 13J30, 90C26; Secondary 12Y05, 13P99, 14P10, 90C22}
\date{\today}

\begin{abstract}
We consider the problem of computing the global infimum of a real
polynomial $f$ on $\R^n$. Every global minimizer of $f$ lies on its gradient
variety, i.e., the algebraic subset of $\R^n$ where the gradient of $f$
vanishes. If $f$ attains a minimum on $\R^n$, it is therefore equivalent to
look for the greatest lower bound of $f$ on its gradient variety. Nie, Demmel
and Sturmfels proved recently a theorem about the existence of sums of squares
certificates for such lower bounds. Based on these certificates, they find
arbitrarily tight relaxations of the original problem that can be formulated
as semidefinite programs and thus be solved efficiently.

We deal here with the more general case when $f$ is bounded from below but
does not necessarily attain a minimum. In this case, the method of Nie, Demmel
and Sturmfels might yield completely wrong results. In order to overcome
this problem, we replace the gradient variety by larger semialgebraic subsets
of $\R^n$ which we call gradient tentacles. It now gets substantially harder
to prove the existence of the necessary sums of squares certificates.
\end{abstract}

\maketitle 

\section{Introduction}

Throughout this article, $\N:=\{1,2,\dots\}$, $\R$ and $\C$ denote the sets
of natural, real and complex numbers, respectively. We fix $n\in\N$,
and consider real polynomials in $n$ variables $\x:=(X_1,\dots,X_n)$. These
polynomials form a commutative ring $$\rx:=\R[X_1,\dots,X_n].$$

\subsection{The problem}

We consider the problem of computing good approximations for the global
infimum $$f^*:=\inf\{f(x)\mid x\in\R^n\}\in\R\cup\{-\infty\}$$
of a polynomial $f\in\rx$. Since $f^*$ is the greatest lower bound of $f$,
it is equivalent to compute
\begin{equation}\label{glb}
f^*=\sup\{a\in\R\mid f-a\ge 0\text{\ on\ }\R^n\}\in\R\cup\{-\infty\}.
\end{equation}
To solve this hard problem, it has become a standard approach to approximate
$f^*$ by exchanging in (\ref{glb}) the nonnegativity constraint
\begin{equation}\label{nnconstraint}
f-a\ge 0\text{\ on\ }\R^n
\end{equation}
by a computationally more feasible condition and analyze the error caused by
this substitution. Typically, the choice of this replacement is related to the
interplay between (globally) nonnegative polynomials, sums of squares of
polynomials and semidefinite optimization (also called semidefinite
programming):

\subsection{Method based on the fact that every sum of squares of polynomials
is nonnegative (Shor \cite{sho}, Stetsyuk \cite{ss}, Parrilo and
Sturmfels \cite{ps} et al.)}\label{shor}

We start with the most basic ideas concerning
these connections which can
be found in greater detail in the just cited references.
A first try is to replace condition (\ref{nnconstraint}) by the constraint
\begin{equation}\label{sosconstraint}
f-a\text{\ is a sum of squares in the polynomial ring $\rx$}
\end{equation}
since every sum of squares in $\rx$ is obviously nonnegative on $\R^n$.

The advantage of (\ref{sosconstraint}) over (\ref{nnconstraint}) is
that sums of squares of polynomials can be nicely parametrized.
Fix a column vector $v$ whose entries are a basis of
the vector space $\rx_d$ of all real polynomials of degree $\le d$
in $n$ variables ($d\in\N_0:=\{0\}\cup\N$). This vector has a certain length
$k=\dim\rx_d$.
It is easy to see that the map from the vector space $S\R^{k\times k}$ of
symmetric $k\times k$-matrices to $\rx_{2d}$ defined by $M\mapsto v^TMv$ is
surjective. Using the spectral theorem for symmetric matrices, it is not hard
to prove that a polynomial $f\in\rx_{2d}$ is a sum of squares in
$\rx$ if and only if $f=v^TMv$ for some positive semidefinite matrix
$M\in S\R^{k\times k}$. Use the following remark which is an easy
exercise (write the polynomials as sums of their homogeneous parts).
\begin{remark}\label{degres}
In any representation $f=\sum_i g_i^2$
of a polynomial $f\in\rx_{2d}$ as a sum of squares $g_i\in\rx$, we have
necessarily $\deg g_i\le d$.
\end{remark}
The described parametrization shows that the modified problem (where we
exchange (\ref{nnconstraint}) by (\ref{sosconstraint})), i.e., the problem
to compute
\begin{equation}\label{sos}
f^\sos:=\sup\{a\in\R\mid f-a\text{ is a sum of squares in $\rx$}\}
\in\R\cup\{-\infty\}
\end{equation}
can be written as a \emph{semidefinite optimization problem}
(also called semidefinite program or SDP for short), i.e., as the problem of
minimizing (or maximizing) an affine
linear function on the intersection of the cone of positive semidefinite
matrices with an affine subspace in $S\R^{k\times k}$. For solving SDPs, there
exist very good numerical algorithms, perhaps almost as good as
for linear optimization problems. Linear optimization can be seen as the
restriction of semidefinite optimization to diagonal matrices, i.e., a method
to minimize an
affine linear function on the intersection of the cone $\nn^k$ with an affine
subspace of $\R^k$. Speaking very vaguely, most concepts from linear
optimization carry over to semidefinite optimization because every symmetric
matrix can be diagonalized. We refer for example to \cite{tod} for an
introduction to semidefinite programming.

Whereas computing $f^*$ as defined in (\ref{glb}) is a very hard problem, it
is relatively easy to compute (numerically to a given precision) $f^\sos$
defined in (\ref{sos}). Of course, the question arises how $f^*$ and
$f^\sos$ are related. Since (\ref{sosconstraint}) implies (\ref{nnconstraint}),
it is clear that $f^\sos\le f^*$. The converse implication (and thus
$f^\sos=f^*$) holds in some cases: A globally nonnegative polynomial
\begin{itemize}
\item in one variable or
\item of degree at most two or
\item in two variables of degree at most four
\end{itemize}
is a sum of squares of polynomials. We refer to \cite{rez} for an overview
of these and related old facts. However,
recently Blekherman has shown in \cite{ble} that for fixed degree $d\ge 4$ and
high number of variables $n$ only a very small portion (in some reasonable
sense) of the globally nonnegative polynomials of degree at most $d$ in $n$
variables are sums of squares. In particular, $f^\sos$ will often differ from
$f^*$. For example, the \emph{Motzkin polynomial}
\begin{equation}\label{motzkin}
M:=X^2Y^2(X^2 + Y^2 - 3Z^2) + Z^6\in\R[X,Y,Z].
\end{equation}
is nonnegative but not a sum of squares (see \cite{rez,ps}). We have
$M^*=0$ but $M^\sos=-\infty$. The latter follows from the fact that $M$ is
homogeneous and not a sum of squares by the following remark applied to
$f:=M-a$ for $a\in\R$ (which can again be proved easily by considering
homogeneous parts).
\begin{remark}\label{homobs}
If $f$ is a sum of squares in $\rx$, then so is the highest homogeneous part
(the leading form) of $f$.
\end{remark}
We see that the basic problem with this method (computing $f^\sos$ by solving
an SDP
and hoping that $f^\sos$ is close to $f^*$) is that polynomials positive on
$\R^n$ in general do not have a representation as a sum of squares, a fact
that Hilbert already knew.

\subsection{The Positivstellensatz}\label{positivstellensatz}

In the 17th of his famous of 23 problems, Hilbert asked whether
every (globally) nonnegative (real) polynomial (in several variables) was
a sum of squares of \emph{rational functions}. Artin answered this
question affirmatively in 1926 and today there exist numerous refinements
of his solution. One of them is the
\emph{Positivstellensatz} (in analogy
to Hilbert's Nullstellensatz). It is often attributed to Stengle \cite{ste}
who clearly deserves credit for finding it independently and making it widely
known. However, Prestel \cite[Section 4.7]{pd} recently discovered that
Krivine \cite{kri} knew the result about ten years earlier in 1964. Here we
state only the following special case of the Positivstellensatz.
\begin{theorem}[Krivine]
For every $f\in\rx$, the following are equivalent.
\begin{enumerate}[(i)]
\item $f>0$ on $\R^n$
\item There are sums of squares $s$ and $t$ in $\rx$ such that $sf=1+t$.
\end{enumerate}
\end{theorem}
By this theorem, we have of course that $f^*$ is the supremum over all $a\in\R$
such that there are sums of squares $s,t\in\rx$ with $s(f-a)=1+t$. When one
tries to write this as an SDP there are two obstacles.

First, each SDP involves matrices of a fixed (finite) size. But with matrices
of a fixed
size, we can only parametrize sums of squares up to a certain degree.
We need therefore to impose a degree restriction on $s$ and $t$. There are no
(at least up to now)
\emph{practically relevant} degree bounds that could guarantee that
such a restriction would not affect the result. We refer to the tremendous
work \cite{scd} of Schmid on degree bounds. 
This first obstacle, namely the question of degrees of the sums of
squares, will us accompany throughout the article. The answer will always be
to model the problem not as a single SDP but as a whole \emph{sequence} of
SDPs, each SDP corresponding to a certain degree restriction. As you solve
one SDP after the other, the degree restriction gets less restrictive and
you hope for fast convergence of the optimal values of the SDPs to $f^*$.
For newcomers in the field, it seems at first glance unsatisfactory
having to deal with a whole sequence of SDPs rather than a single SDP. But
after all, it is only natural that a very hard problem cannot be modeled by an
SDP of a reasonable size so that you have to look for good
\emph{relaxations} of the problem which can easier be dealt with and to
which the techniques of mathematical optimization can be applied.

The second obstacle is much more severe. It is the fact that the unknown
polynomial $s\in\rx$ is multiplied with the unknown $a\in\R$ on the left hand
side of the constraint $s(f-a)=1+t$. This makes the formulation as an SDP
(even after having imposed a restriction on the degree of $s$ and $t$)
impossible
(or at least highly non-obvious). Of course, if you \emph{fix} $a\in\R$ and
a degree bound $2d$ for $s$ and $t$, then the question whether there exist
sums of squares $s$ and $t$ of degree at most $2d$ such that $s(f-a)=1+t$
is equivalent to the feasibility of an SDP. But this plays (at least currently)
only a role as a criterion that \emph{might help to decide} whether a
certain fixed (or guessed) $a\in\R$ is a strict lower bound of $f$. We refer to
\cite{ps} for more details. What one needs are representation theorems for
positive polynomials that are better suited for optimization than the
Positivstellensatz (even if they are sometimes less aesthetic).

\subsection{``Big ball'' method proposed by Lasserre \cite{l1}}\label{bigball}

In the last 15 years, a lot of progress has been made in proving existence
of sums of squares certificates which can be exploited for optimization
(although most of the new results were obtained without having in mind the
application in optimization which has been established more recently).
The first breakthrough was perhaps Schm\"udgen's theorem
\cite[Corollary 3]{sch} all of whose proofs use the Positivstellensatz.
In this article, we will prove a generalization of Schm\"udgen's theorem,
namely Theorem \ref{schmuedgen} below. In \cite{l1}, Lasserre uses the
following special case of Schm\"udgen's theorem which has already been proved
by Cassier \cite[Th\'eor\`eme 4]{cas} and which can even be derived easily
from \cite[Th\'eor\`eme 12]{kri}.
\begin{theorem}[Cassier]\label{cassier}
For $f\in\rx$ and $R\ge 0$, the following are equivalent.
\begin{enumerate}[(i)]
\item $f\ge 0$ on the closed ball centered at the origin of radius $R$
\item For all $\ep>0$, there are sums of squares $s$ and $t$ in $\rx$ such that
$$f+\ep=s+t(R^2-\|\x\|^2).$$
\end{enumerate}
\end{theorem}
Here and in the following, we use the notation
$$\|\x\|^2:=X_1^2+\dots+X_n^2\in\rx.$$
Similar to Subsection \ref{shor}, it can be seen that for any fixed
$d\in\N_0$,
computing the supremum over all $a\in\R$ such that $f-a=s+t(R^2-\|\x\|^2)$ for
some sums of squares $s,t\in\rx$ of degree at most $2d$ amounts to solving an
SDP. Therefore you get a sequence of SDPs parametrized by $d\in\N_0$.
Theorem \ref{cassier} can now be interpreted as a convergence result, namely
the sequence of optimal values of these SDPs converges to the minimum of $f$
on the closed ball around the origin with radius $R$. If one has a polynomial
$f\in\rx$ attaining a minimum on $\R^n$ and for which one knows moreover a big
ball on which this minimum is attained, this method is good for computing
$f^*$. Of course, if you do not know such a big ball in advance you might
choose larger and larger $R$. But at the same time you might have to choose
a bigger and bigger degree restriction $d\in\N_0$ and it is not really clear
how to get a sequence of SDPs that converges to $f^*$.

\subsection{Lasserre's high order perturbation method \cite{l2}}

Recently, Lasserre used in \cite{l2} a theorem of Nussbaum from operator
theory to prove the following result that can be exploited in a similar way
for global optimization of polynomials.
\begin{theorem}[Lasserre]\label{lasserre}
For every $f\in\rx$, the following are equivalent:
\begin{enumerate}[(i)]
\item \label{nnrn}
      $f\ge 0$ on $\R^n$
\item \label{expappr}
      For all $\ep>0$, there is $r\in\N_0$ such that
$$f+\ep\sum_{i=1}^n\sum_{k=0}^r\frac{X_i^{2k}}{k!}\text{\ is a sum of squares
in $\rx$}.$$
\end{enumerate}
\end{theorem}
Note that (\ref{expappr}) implies that $f(x)+\ep\sum_{i=1}^n\exp(x_i)\ge 0$
for all $x\in\R^n$ and $\ep>0$ which in turn implies (\ref{nnrn}). In condition
(\ref{expappr}), $r$ depends on $\ep$ and $f$. Using real algebra and model
theory, Netzer showed that in fact $r$ depends only on $\ep$, $n$, the
degree of $f$ and a bound on the size of the coefficients of $f$
\cite{net,ln}.

\subsection{``Gradient perturbation'' method proposed by Jibetean and Laurent
            \cite{jl}}\label{jlsubsection}

The most standard idea for finding the minimum of a function everybody knows
from calculus is to compute critical points, i.e., the points where the
gradient vanishes. It is a natural question
whether the power of classical differential calculus can be combined with
the relatively new ideas using sums of squares. Fortunately, it can and the
rest of the article will be about how to merge both concepts, sums of squares
and differential calculus.

If a polynomial $f\in\rx$ attains a minimum in $x\in\R^n$, i.e.,
$f(x)\le f(y)$ for all $y\in\R^n$, then the gradient $\nabla f$
of $f$ vanishes at $x$, i.e., $\nabla f(x)=0$. However, there are polynomials
that are bounded from below on $\R^n$ and yet do not attain a minimum on
$\R^n$. The simplest example is perhaps
\begin{equation}\label{doesnot}
f:=(1-XY)^2+Y^2\in\R[X,Y]
\end{equation}
for which we have $f>0$ on $\R^n$ but $f^*=0$ since
$\lim_{x\to\infty}f(x,\frac 1x)=0$.  In the following,
$$(\nabla f):=\left(\frac{\partial f}{\partial X_1},\dots,
  \frac{\partial f}{\partial X_n}\right)\subseteq\rx$$
denotes the ideal generated by the partial derivatives of $f$ in $\rx$.
We call this ideal the \emph{gradient ideal} of $f$.

Without going into details, the basic idea of
Jibetean and Laurent in \cite{jl} is again to apply a perturbation to $f$.
Instead of adding a truncated exponential like Lasserre, they just add
$\ep\sum_{i=1}^nX_i^{2(d+1)}$ for small $\ep>0$ when $\deg f=2d$.
If $f>0$ on $\R^n$, then the perturbed polynomial
$f_\ep:=f+\ep\|\x\|^{2(d+1)}$ is again a sum of squares
but this time only modulo its gradient ideal $(\nabla f_\ep)$.
In this case, this is quite easy to prove since it turns out that 
this ideal will be zero-dimensional, i.e.,
$\rx/(\nabla f_\ep)$ is a finite-dimensional real algebra.
We will later see in Theorems \ref{nie} and \ref{zeltschmid} that this
finite-dimensionality is not needed for the sums of squares representation.
But the work of Jibetean and Laurent exploits the finite-dimensionality in
many ways. We refer to \cite{jl} for details.

\subsection{``Gradient variety'' method by Nie, Demmel and Sturmfels
            \cite{nds}}\label{gradvariety}

The two perturbation methods just sketched rely on introducing very small
coefficients in a polynomial. This small coefficients might lead to SDPs
which are hard to solve because of numerical instability. It is therefore
natural to think of another method which avoids perturbation at all.
Nie, Demmels and Sturmfels considered, for a polynomial $f\in\rx$, its
\emph{gradient variety}
$$V(\nabla f):=\{x\in\C^n\mid\nabla f(x)=0\}.$$
This is the algebraic variety corresponding to the radical of the
gradient ideal $(\nabla f)$.
It can be shown that a polynomial $f\in\rx$ is constant on each irreducible
component of the gradient variety (see \cite{nds} or use an unpublished
algebraic argument of Scheiderer based on K\"ahler differentials). This is
the key to show that a polynomial $f\in\rx$ nonnegative on its gradient
variety is a sum of squares modulo its gradient ideal in the case where the
ideal is radical. In the general case where the gradient ideal is not
necessarily radical, the same thing still holds for polynomials
\emph{positive} on their gradient variety. The following is
essentially \cite[Theorem 9]{nds} (confer also the recent work \cite{m2}).
We will later prove a generalization of this theorem as a byproduct. See
Corollary \ref{gradpos} below.

\begin{theorem}[Nie, Demmel and Sturmfels]\label{nie}
For every $f\in\rx$ attaining a minimum on $\R^n$, the following are
equivalent.
\begin{enumerate}[(i)]
\item \label{niern}
$f\ge 0$ on $\R^n$
\item \label{nienabla}
$f\ge 0$ on $V(\nabla f)\cap\R^n$
\item \label{nierep}
For all $\ep>0$, there exists a sum of squares $s$ in $\rx$
such that $$f+\ep\in s+(\nabla f).$$
\end{enumerate}
Moreover, (\ref{nienabla}) and (\ref{nierep}) are equivalent for all $f\in\rx$.
\end{theorem}
For each degree restriction $d\in\N_0$, the problem of computing the supremum
over all $a\in\R$ such that
$$f-a=s+p_1\frac{\partial f}{\partial X_1}+\dots+
  p_n\frac{\partial f}{\partial X_n}$$
for some sum of squares $s$ in $\rx$ and polynomials $p_1,\dots,p_n$ of degree
at most $d$, can be expressed as an SDP. Theorem \ref{nie} shows that the
optimal values of the corresponding sequence of SDPs (indexed by $d$)
tend to $f^*$ provided that $f$ attains a minimum on $\R^n$. However, if $f$
does not attain a minimum on $\R^n$, the computed sequence still tends to the
infimum of $f$ on its gradient variety which might however now be very
different from $f^*$. Take for example the polynomial $f$
from (\ref{doesnot}). It is easy to see that $V(\nabla f)=\{0\}$ and therefore
the method computes $f(0)=1$ instead of $f^*=0$.
In \cite[Section 7]{nds}, the authors write:
\begin{quotation}
``This paper proposes a method for minimizing a multivariate polynomial
$f(x)$ over its gradient variety. We assume that the infimum $f^*$
is attained. This assumption is non-trivial, and we do not address the
(important and difficult) question of how to verify that a given polynomial
$f(x)$ has this property.''
\end{quotation}

\subsection{Our ``gradient tentacle'' method}

The reason why the method just described might fail is that the global
infimum of a polynomial $f\in\rx$ is not always a \emph{critical value}
of $f$, i.e., a value that $f$ takes on at least one of its critical points in
$\R^n$. Now there is a well-established notion of
\emph{generalized critical values} which includes also the
\emph{asymptotic critical values} (a kind of critical values at
infinity we will introduce in Definition \ref{critvaldef} below).

In this article, we will replace the real part $V(\nabla f)\cap\R^n$
of the gradient variety by several larger semialgebraic sets on which
the partial derivatives do not necessarily vanish but get very small
far away from the origin. These semialgebraic sets often look
like tentacles, and that is how we will call them. All tentacles we will
consider are defined by a single polynomial inequality that depends only
on the polynomial
$$\|\nabla f\|^2:=\left(\frac{\partial f}{\partial X_1}\right)^2+\dots+
                  \left(\frac{\partial f}{\partial X_n}\right)^2$$
and expresses that this polynomial gets very small. Given a polynomial $f$ for
which you want to compute $f^*$, the game will consist in finding a tentacle
such that two things will hold at the same time:
\begin{itemize}
\item There exist suitable sums of squares certificates for nonnegativity
on the tentacle.
\item The infimum of $f$ on $\R^n$ and on the tentacle coincide.
\end{itemize}
One can imagine that these two properties are hardly compatible. Taking
$\R^n$ as a tentacle, would of course ensure the second condition but we
have discussed in Subsection \ref{shor} that the first one would be badly
violated. The other extreme would be to take the empty set as a tentacle.
Then the first condition would trivially be satisfied whereas the second
would fail badly. How we will roughly be able to find the balancing act
between the two requirements is as follows: The second condition will be
satisfied by known non-trivial theorems about asymptotic behaviour of
polynomials at infinity. The existence of suitable sums of squares
certificates will be based on the author's (real) algebraic work \cite{sr1} on
iterated rings of bounded elements (also called real holomorphy rings).

\subsection{Contents of the article} The article is organized as follows.
In Section \ref{sossection}, we prove a general sums of squares
representation theorem which generalizes Schm\"udgen's theorem we have
mentioned in Subsection \ref{bigball}. This representation theorem is
interesting in itself and will be used in the subsequent sections.
In Section \ref{parusinskisection}, we introduce
a gradient tentacle (see Definition \ref{ptentdef}) which is defined by the
polynomial inequality
$$\|\nabla f\|^2\|\x\|^2\le 1.$$
We call this gradient tentacle \emph{principal} since we can prove that
it does the job in a large number of cases (see Theorem \ref{ptheorem}) and
there is hope that it works in fact for all polynomials $f\in\rx$
bounded from below. Indeed, we have not found any counterexamples
(see Open Problem \ref{open}). In case this hope were disappointed, we
present in Section \ref{highersection}
a collection of other gradient tentacles (see
Definition \ref{tentdef}) defined by the polynomial inequalities
$$\|\nabla f\|^{2N}(1+\|\x\|^2)^{N+1}\le 1\qquad (N\in\N).$$
Their advantage is that if $f\in\rx$ is bounded from below and $N$ is large
enough for this particular $f$, then we can prove that the corresponding
tentacle does the job (see Theorems \ref{zeltschmid} and \ref{opt}). We call
these tentacles
higher gradient tentacles since the degree of the defining inequality gets
unfortunately high when $N$ gets big which has certainly negative consequences
for the complexity of solving the SDPs arising from these tentacles.
However, if $f$ attains a minimum on $\R^n$, then any choice of $N\in\N$ will
be good. Conclusions are drawn in Section \ref{conclusion}.

\section{The sums of squares representation}\label{sossection}

In this section, we prove the important sums of squares representation theorem
we will need in the following sections. It is a generalization of
Schm\"udgen's Positivstellensatz (see \cite{pd,sch}) which is also of
independent interest. Schm\"udgen's result is not to confuse
with the (classical) Positivstellensatz we described in the introduction.
The connection between the two is that all known proofs of Schm\"udgen's
result use the classical Positivstellensatz. Our result, Theorem
\ref{schmuedgen} below, is much harder to prove than Schm\"udgen's result.
Its proof relies on the theory of iterated
\emph{rings of bounded elements}
(also called real holomorphy rings) described in \cite{sr1}.

\begin{definition}
For any polynomial $f\in\rx$ and subset $S\subseteq\R^n$, the set
$\ran(f, S)$ of \emph{asymptotic values} of
$f$ on $S$ consists of all $y\in\R$ for which there exists a sequence
$(x_k)_{k\in\N}$ of points $x_k\in S$ such that
\begin{equation}\label{asymptotic}
\lim_{k\to\infty}\|x_k\|=\infty\qquad\text{and}\qquad
\lim_{k\to\infty}f(x_k)=y.
\end{equation}
\end{definition}

We now recall the important notion of a preordering of a commutative ring.
Except in the proof of Theorem \ref{schmuedgen}, we need this concept
only for the ring $\rx$.

\begin{definition}
Let $A$ be a commutative ring (with $1$).
A subset $T\subseteq A$ is called a \emph{preordering} if it contains
all squares $f^2$ of elements $f\in A$ and is closed under addition and
multiplication. The preordering \emph{generated} by $g_1,\dots,g_m\in A$
\begin{equation}\label{explicit}
T(g_1,\dots,g_m)=\left\{\sum_{\de\in\{0,1\}^m}s_\de g_1^{\de_1}\dots
g_m^{\de_m}\mid s_\de\text{\ is a sum of squares in\ } A\right\}
\end{equation}
is by definition the smallest preordering containing $g_1,\dots,g_m$.
\end{definition}

If $g_1,\dots,g_m\in\rx$ are polynomials, then
the elements of $T(g_1,\dots,g_m)$ have obviously the geometric
property that they are nonnegative on the
(basic closed semialgebraic) set $S$ they define by (\ref{basicclosed}) below.
The next theorem is a partial converse. Namely, if a polynomial satisfies
on $S$ some stronger geometric condition, then it lies necessarily
in $T(g_1,\dots,g_m)$. In case that $S$ is compact, the conditions
(\ref{sbounded}) and (\ref{sfinite}) below are empty and the theorem is
Schm\"udgen's Positivstellensatz (see \cite{pd,sch}). The more general
version we need here is quite hard to prove.

\begin{theorem}\label{schmuedgen}
Let $f,g_1,\dots,g_m\in\rx$ and set
\begin{equation}\label{basicclosed}
S:=\{x\in\R^n\mid g_1(x)\ge 0,\dots,g_m(x)\ge 0\}.
\end{equation}
Suppose that
\begin{enumerate}[(a)]
\item $f$ is bounded on $S$,\label{sbounded}
\item $f$ has only finitely many asymptotic values on $S$ and all of
these are positive, i.e., $\ran(f,S)$ is a finite subset of $\R_{>0}$, and
\label{sfinite}
\item $f>0$ on $S$.\label{spositive}
\end{enumerate}
Then $f\in T(g_1,\dots,g_m)$.
\end{theorem}

\begin{proof}
Write $\ran(f,S)=\{y_1,\dots,y_s\}\subseteq\R_{>0}$ and consider the
polynomial $$h:=\prod_{i=1}^s(f-y_i).$$
This polynomial is ``on $S$ small at infinity'' by which we mean that
for every $\ep>0$ there exists $k\in\N$ such that for all $x\in S$
with $\|x\|\ge k$, we have $|h(x)|<\ep$.

To show this, assume the contrary. Then there exists $\ep>0$ and a sequence
$(x_k)_{k\in\N}$ of points $x_k\in S$ with
$\lim_{k\to\infty}\|x_k\|=\infty$ and
\begin{equation}\label{stayaway}
|h(x_k)|\ge\ep\qquad\text{for all $k\in\N$.}
\end{equation}
Because the sequence $(f(x_k))_{k\in\N}$ is bounded by hypothesis
(\ref{sbounded}), we find an infinite subset $I\subseteq\N$
such that the subsequence $(f(x_k))_{k\in I}$ converges. The limit must
be one of the asymptotic values of $f$ on $S$, i.e.,
$\lim_{k\in I,k\to\infty}f(x_k)=y_i$ for some $i\in\{1,\dots,s\}$. Using
(\ref{sbounded}), it follows that
$\lim_{k\in I,k\to\infty}h(x_k)=0$ contradicting (\ref{stayaway}).

Let $A:=(\rx,T)$ where $T:=T(g_1,\dots,g_m)$. The set
$$H'(A):=\{p\in\rx\mid N\pm p\in T\text{\ for some $N\in\N$}\}$$
is a subring of $A$ (see, e.g, \cite[Definition 1.2]{sr1}).
We endow $H'(A)$ with the preordering $T':=T\cap H'(A)$ and consider it as
also as a preordered ring. By \cite[Corollary 3.7]{sr1}, the smallness
of $h$ at infinity proved above is equivalent to $h\in S_\infty(A)$ in the
notation of \cite{sr1}. By \cite[Corollary 4.17]{sr1}, we have
$S_\infty(A)\subseteq H'(A)$ and consequently $h\in H'(A)$. The advantage of
$H'(A)$ over $A$ is that its preordering is archimedean, i.e., $T'+\Z=H'(A)$.
According to an old criterion for an element to be contained in an archimedean
preordering (see for example \cite[Proposition 5.2.3 and Lemma 5.2.7]{pd} or
\cite[Theorem 1.3]{sr1}), our claim $f\in T'$ follows if we can show that
$\ph(f)>0$ for all ring homomorphisms $\ph:H'(A)\to\R$ with
$\ph(T')\subseteq\nn$. For all such homomorphisms possessing an
extension $\bar\ph:A\to\R$ with $\bar\ph(T)\subseteq\nn$, this follows
from hypothesis (\ref{spositive}) because it is easy to see that such an
extension $\bar\ph$ must be evaluation $p\mapsto p(x)$ in the point
$x:=(\bar\ph(X_1),\dots,\bar\ph(X_n))\in S$. Using the theory in \cite{sr1},
we will see that the only possibility for such a $\ph$ not to have such an
extension $\bar\ph$ is that $\ph(h)=0$. Then we will be done since
$\ph(h)=0$ implies $\ph(f)=y_i>0$ for some $i$. We have used here that
$f\in H'(A)$ which follows from $h\in H'(A)$ since $H'(A)$ is integrally
closed in $A$ (see \cite[Theorem 5.3]{sr1}).

So let us now use \cite{sr1}. By \cite[Corollary 3.7 and Theorem 4.18]{sr1},
the smallness of $h$ at infinity means that
$$A_h=H'(A)_h$$
where we deal on both sides of this equation with the localization of a
preordered ring by the element $h$ (see \cite[pages 24 and 25]{sr1}).
If $\ph:H'(A)\to\R$ is a ring homomorphism with $\ph(T')\subseteq\nn$
and $\ph(h)\neq 0$, then $\ph$ extends to a ring homomorphism
$\tilde\ph:A_h=H'(A)_h\to\R$ with
$\tilde\ph(T_h)=\tilde\ph(T'_h)\subseteq\nn$.
Then $\bar\ph:=\tilde\ph|_A$ is the desired extension of $\ph$.
\end{proof}

\begin{example}\label{hnex}
Consider the polynomials
\begin{equation}\label{hndef}
h_N:=1-Y^N(1+X)^{N+1}\in\R[X,Y]\qquad(N\in\N)
\end{equation}
in two variables. We fix $N\in\N$ and apply Theorem \ref{schmuedgen} with
$f=h_{N+1}$, $m=3$, $g_1=X$, $g_2=Y$ und $g_3=h_N$. The set $S$ defined by
the $g_i$ as in (\ref{basicclosed}) is a subset of the first quadrant which
is bounded in $Y$-direction but unbounded in $X$-direction. Of course, we have
$0\le h_N\le 1$ and
$$0\le Y(1+X)\le\frac 1{\sqrt[N]{1+X}}\qquad\text{on $S$}$$ showing that $0$
is the only asymptotic
value of $$1-h_{N+1}=(1-h_N)Y(1+X)$$ on $S$
and therefore $\ran(h_{N+1},S)=\{1\}$. It follows also that
$0\le h_{N+1}\le 1$ on $S$. By Theorem \ref{schmuedgen}, we obtain
\begin{equation}\label{hnweakrep}
h_{N+1}+\ep\in T(X,Y,h_N)
\end{equation}
for all $\ep>0$.
\end{example}

The following lemma shows that (\ref{hnweakrep}) holds even for $\ep=0$,
a fact that does not follow from Theorem \ref{schmuedgen}. This lemma will
be interesting later to compare the quality of certain SDP relaxations
(see Proposition \ref{tentinccert}). In its proof, we will explicitly construct
a representation of $h_{N+1}$ as an element of $T(X,Y,h_N)$. Only part of this
explicit representation will be needed in the sequel, namely an explicit
polynomial $g\in T(X,Y)$ such that $h_{N+1}\in T(X,Y)+gh_N\subseteq
T(X,Y,h_N)$. This explains the formulation of the statement.
Theorem \ref{schmuedgen} will not be used in the proof but gave us good hope
before we had the proof. The role of Theorem \ref{schmuedgen} in this article
is above all to prove Theorems \ref{ptheorem} and \ref{zeltschmid} below. 

\begin{lemma}\label{hnlemma}
For the polynomials $h_N$ defined by (\ref{hndef}), we have
$$h_{N+1}-\left(1+\frac 1N\right)Y(1+X)h_N\in T(X,Y).$$
\end{lemma}

\begin{proof}
For a new variable $Z$,
\begin{align*}
(Z-1)^2\sum_{k=0}^{N-1}(N-k)Z^k&=(Z-1)^2\left(N\sum_{k=0}^{N-1}Z^k-Z
\sum_{k=1}^{N-1}kZ^{k-1}\right)\\
&=(Z-1)^2\left(N\frac{Z^N-1}{Z-1}-Z\frac\partial{\partial Z}
\left(\frac{Z^N-1}{Z-1}\right)\right)\\
&=N(Z-1)(Z^N-1)-Z((Z-1)NZ^{N-1}-(Z^N-1))\\
&=Z^{N+1}-(N+1)Z+N.
\end{align*}
Specializing $Z$ to $z:=Y(1+X)$, we have therefore
\begin{align*}
Nh_{N+1}-(N+1)zh_N
&=N(1-z^{N+1}(1+X))-(N+1)z(1-z^N(1+X))\\
&=z^{N+1}X+(z^{N+1}-(N+1)z+N)\\
&=z^{N+1}X+(z-1)^2\sum_{k=0}^{N-1}(N-k)z^k\in T(X,Y).
\end{align*}
Dividing by $N=(\sqrt N)^2$ yields our claim.
\end{proof}

\section{The principal gradient tentacle}\label{parusinskisection}

In this section, we associate to every polynomial $f\in\rx$ a gradient
tentacle which is a subset of $\R^n$ containing the real part of the gradient
variety of $f$ and defined by a single polynomial inequality
whose degree is not more than twice the degree of $f$. The infimum of any
polynomial $f\in\rx$ bounded from below on $\R^n$ will coincide with the
infimum on its principal gradient tentacle (see Theorem \ref{critical}).
Under some technical assumption (see Definition \ref{isolated}) which is not
known to be necessary (see Open Problem \ref{open}), we prove a sums of squares
certificate for nonnegativity of $f$ on its principal gradient tentacle which
is suitable for optimization purposes. This representation theorem (Theorem
\ref{ptheorem}) is of independent interest and its proof is mainly based
on the nontrivial representation theorem from the previous section and
a result of Parusi\'nski on the behaviour of polynomials at infinity
(\cite[Theorem 1.4]{p1}). In Subsection \ref{psdp}, we outline how to get
a sequence of SDPs growing in size whose optimal values tend to $f^*$ for any
$f$ satisfying the conditions of Theorem \ref{ptheorem} (or perhaps for any
$f$ with
$f^*>-\infty$ if the answer to Open Problem \ref{open} is yes). In the
Subsections \ref{yalmip} and \ref{sostools}, we give a MATLAB code for the
sums of squares optimization toolboxes YALMIP \cite{löf} and SOSTOOLS
\cite{pps} that produces and solves these SDP relaxations. This short and
simple code is meant for readers who have little experience with such
toolboxes and want nevertheless try our proposed method on their own. In
Subsection \ref{numerical}, we provide simple examples which have been
calculated using the YALMIP code from Subsection \ref{yalmip}.

We start by recalling the concept of asymptotic critical values developed by
Rabier in his 1997 milestone paper \cite{rab}. For simplicity, we stay in the
setting of real polynomials right from the beginning
(though part of this theory make sense in a much broader context).

\begin{definition}\label{critvaldef}
Suppose $f\in\rx$. The set $K_0(f)$ of \emph{critical values}
of $f$ consists of all $y\in\R$ for which there exists $x\in\R^n$ such that
$\nabla f(x)=0$ and $f(x)=y$. The set $K(f)$ of
\emph{generalized critical values}
of $f$ consists of all $y\in\R$ for which there exists a sequence
$(x_k)_{k\in\N}$ in $\R^n$ such that
\begin{equation}\label{gradr}
\lim_{k\to\infty}\|\nabla f(x_k)\|(1+\|x_k\|)=0\qquad\text{and}\qquad
\lim_{k\to\infty}f(x_k)=y.
\end{equation}
The set $K_\infty(f)$ of \emph{asymptotic critical values}
consists of all
$y\in\R$ for which there exists a sequence $(x_k)_{k\in\N}$ in $\R^n$ such
that $\lim_{k\to\infty}\|x_k\|=\infty$ and (\ref{gradr}) hold.
\end{definition}

The following proposition is easy.

\begin{proposition}\label{union}
The set of generalized critical values of a polynomial
$f\in\rx$ is the union of its set of critical and asymptotic
critical values, i.e., $$K(f)=K_0(f)\cup K_\infty(f).$$
\end{proposition}

The following notions go back to Thom \cite{tho}.

\begin{definition}\label{typical}
Suppose $f\in\rx$. We say that $y\in\R$ is a \emph{typical value} of $f$ if
there is neighbourhood $U$ of $y$ in $\R$ and a smooth
(i.e., $C^\infty$) manifold $F$ such that $f|_{f^{-1}(U)}:f^{-1}(U)\to U$ is a
(not necessarily surjective) trivial smooth fiber bundle, i.e., there exists a
smooth manifold $F$ and a
$C^\infty$ diffeomorphism $\Ph:f^{-1}(U)\to F\times U$ such that
$f|_{f^{-1}(U)}=\pi_2\circ\Ph$ where $\pi_2:F\times U\to U$ is the canonical
projection. We call $y\in\R$ an \emph{atypical value} of $f$ if it is not a
typical value of $f$. The set of all atypical values of $f$ is denoted by
$B(f)$ and called the \emph{bifurcation set} of $f$.
\end{definition}

Note that a $\Ph$ like in the above definition induces a $C^\infty$
diffeomorphism $f^{-1}(y)\to F\times\{y\}\cong F$ for every $y\in U$.
In this context, the preimages $f^{-1}(y)$ are called fibers and $F$ is
called \emph{the} fiber. We do not require that the fiber bundle
$f|_{f^{-1}(U)}:f^{-1}(U)\to U$ is surjective (if it is not then the image is
necessarily empty). Hence the fiber $F$ may be empty
and a \emph{typical value} is not necessarily a value taken on by $f$.
We make use the following well-known theorem
(see, e.g., \cite[Theorem 3.1]{kos}).

\begin{theorem}\label{bifurcation}
Suppose $f\in\rx$. Then $B(f)\subseteq K(f)$ and $K(f)$ is finite.
\end{theorem}

The advantage of $K(f)$ over $K_0(f)$ is that $f^*\in K(f)$ even if $f$ does
not attain a minimum on $\R^n$. This is an easy consequence of Theorem
\ref{bifurcation}. See Theorem \ref{critical} below.

\begin{example}\label{bex}
Consider again the polynomial $f=(1-XY)^2+Y^2\in R[X,Y]$ from
(\ref{doesnot}) that does not attain its infimum $f^*=0$ on $\R^2$.
Calculating the partial derivatives, it is easy to see that the origin is the
only critical
point of $f$. Because $f$ takes the value $1$ at the origin, we have
$K_0(f)=\{1\}$ and therefore $f^*=0\notin K_0(f)$. Clearly, we have $0\in B(f)$
since $f^{-1}(-y)=\emptyset\neq f^{-1}(y)$ for small $y\in\R_{>0}$.
By Theorem \ref{bifurcation}, we have therefore $0\in K_\infty(f)\subseteq
K(f)$. To show this directly, a first guess would be that
$\|\nabla f(x,\frac 1x)\|(1+\|(x,\frac 1x)\|)$ tends to zero when $x\to\infty$
because $\lim_{x\to\infty}f(x,\frac 1x)=0$. But in fact, this expressions
tends to
$2$ when $x\to\infty$. However, a calculation shows that
$\lim_{x\to\infty}\|\nabla f(x,\frac 1x)\|(1+\|(x,\frac 1x-\frac 1{x^3})\|)=0$.
\end{example}

\begin{definition}\label{ptentdef}
For a polynomial $f\in\rx$, we call
$$\ptent f:=\{x\in\R^n\mid\|\nabla f(x)\|\|x\|\le 1\}$$
the \emph{principal gradient tentacle} of $f$.
\end{definition}

\begin{remark}\label{nierem}
In the definition of $\ptent f$, the inequality $\|\nabla f(x)\|\|x\|\le 1$
could be exchanged by $\|\nabla f(x)\|\|x\|\le R$ for some constant $R>0$.
Then all subsequent results will still hold with obvious modifications.
Using an $R$ different from $1$ might have in certain cases a practical
advantage (see Subsection \ref{stability} below). However, we decided to stay
with this definition in order to get not too technical and to keep the paper
readable.
\end{remark}

As expressed by the notation $\ptent f$, polynomials $f$ with the same
gradient $\nabla f$ have the same gradient tentacle, in other words
$$\ptent{(f+a)}=\ptent f\qquad\text{for all $a\in\R$.}$$

The first important property of $\ptent f$ is stated in the following
immediate consequence of Theorem \ref{bifurcation}.

\begin{theorem}\label{critical}
Suppose $f\in\rx$ is bounded from below. Then $f^*\in K(f)$ and
therefore $f^*=\inf\{f(x)\mid x\in\ptent f\}$.
\end{theorem}

\begin{proof}
By Theorem \ref{bifurcation}, it suffices to show that $f^*\in B(f)$.
Assume that $f^*\notin B(f)$, i.e., $f^*$ is a typical value of $f$. Then
for all $y$ in a neighbourhood of $f^*$, the fibers $f^{-1}(y)$ are smoothly
diffeomorphic to each other. But this is absurd since $f^{-1}(y)$ is empty
for $y<f^*$ but certainly not empty in a neighbourhood of $f^*$.
\end{proof}

Let $\mathbb P^{n-1}(\C)$ denote the $(n-1)$-dimensional complex projective
space over $\C$. For a homogeneous polynomial $f$ and a point
$z\in\mathbb P^{n-1}(\C)$, we simply say $f(z)=0$ to express that $f$ vanishes
on (a non-zero point of) the straight line $z\subseteq\C^n$.
Following \cite{p1}, we give the following definition.

\begin{definition}\label{isolated}
We say that a polynomial $f\in\C[\x]$ \emph{has only isolated singularities at
infinity} if $f\in\C$ (i.e., $f$ is constant) or $d:=\deg f\ge 1$ and
there are only finitely many $z\in\mathbb P^{n-1}(\C)$ such that
\begin{equation}\label{singularity}
\frac{\partial f_d}{\partial X_1}(z)=\dots=
\frac{\partial f_d}{\partial X_n}(z)=f_{d-1}(z)=0
\end{equation}
where $f=\sum_i f_i$ and each $f_i\in\C[\x]$ is zero or homogeneous of
degree $i$.
\end{definition}

As shown in \cite[Section 1.1]{p1}, the geometric interpretation of the above
definition is that the projective closure of a generic fiber of $f$ has only
isolated singularities.

\begin{remark}\label{generic}
A generic complex polynomial has only isolated singularities at infinity. In
fact, much
more is true: A generic polynomial $f\in\C[\x]$ of degree $d\ge 1$ has
\emph{no} isolated singularities at infinity in the sense that there is no
$z\in\mathbb P^{n-1}(\C)$ such that \eqref{singularity} holds.
In more precise words, to every $d\ge 2$, there exists a complex
polynomial relation that is valid for all coefficient tuples of
polynomials $f\in\C[\x]$ of degree $d$ for which \eqref{singularity}
has an infinite number of solutions. This follows from the fact that for a
generic homogeneous polynomial $g\in\C[\x]$ of degree $d\ge 1$, there are only
finitely many points $z\in\mathbb P^{n-1}(\C)$ such that
$\frac{\partial f}{\partial X_i}(z)=0$ for all $i$. See
\cite[Th\'eor\`eme II]{kus} or \cite[Proposition 1.1.1]{shu}.
\end{remark}

\begin{remark}\label{gen2}
In the two variable case $n=2$, every polynomial $f\in\C[\x]$ has only isolated
singularities at infinity. This is clear since \eqref{singularity} defines
an algebraic subvariety of $\mathbb P^1(\C)$.
\end{remark}

The following theorem follows easily from \cite[Theorem 1.4]{p1}.

\begin{theorem}\label{parusinski}
Suppose $f\in\rx$ has only isolated singularities at infinity. Then
$$\ran(f,\ptent f)\subseteq K(f).$$
In particular, $\ran(f,\ptent f)$  is finite, i.e.,
$f$ has only finitely many asymptotic values on its
principal gradient tentacle.
\end{theorem}

\begin{proof}
Let $(x_k)_{k\in\N}$ be a sequence of points $x_k\in\ptent f$ and $y\in\R$
such that $\lim_{k\to\infty}||x_k||=\infty$ and
$\lim_{k\to\infty}f(x_k)=y\notin K_0(f)$. We show that $y\in K_\infty(f)$
using implication (i)$\implies$(ii) in \cite[Theorem 1.4]{p1}.
Because of our sequence $(x_k)_{k\in\N}$, it is impossible that there exists
$N\ge 1$ and $\de >0$ such that for all $x\in\R^n$ with $\|x\|$ sufficiently
large and $f(x)$ sufficiently close to $y$, we have
$$||x||||\nabla f(x)\|\ge\de\sqrt[n]{\|x\|}.$$
This means that condition (ii) in \cite[Theorem 1.4]{p1} is violated.
The implication (i)$\implies$(ii) in \cite[Theorem 1.4]{p1} yields that
$y\in B(f)$ (here we use that $y\notin K_0(f)$). But $B(f)\subseteq K(f)$ by
Theorem \ref{bifurcation}. This shows
$y\in K(f)\setminus K_0(f)\subseteq K_\infty(f)$
by Proposition \ref{union}.
\end{proof}

\begin{lemma}\label{lojasiewicz}
Every $f\in\rx$ is bounded on $\ptent f$.
\end{lemma}

\begin{proof}
By the \L ojasiewicz inequality at infinity \cite[Theorem 1]{spo}, there exist
$c_1,c_2\in\N$ such that for all $x\in\C^n$,
$$|f(x)|\ge c_1\implies|f(x)|\le c_2\|\nabla f(x)\|\|x\|.$$
Then $|f|\le\max\{c_1,c_2\}$ on $\ptent f$.
\end{proof}

\subsection{The principal gradient tentacle and sums of squares}

Here comes one of the main results of this article which is interesting on
its own but can later be read as a convergence result for a sequence of optimal
values of SDPs (Theorem \ref{popt} below).

\begin{theorem}\label{ptheorem}
Let $f\in\rx$ be bounded from below. Furthermore, suppose that $f$ has only
isolated singularities at infinity (which is always true in the two variable
case $n=2$) or the principal gradient tentacle $\ptent f$ is compact.
Then the following are equivalent.
\begin{enumerate}[(i)]
\item $f\ge 0$ on $\R^n$\label{pglobnn}
\item $f\ge 0$ on $\ptent f$\label{ptentnn}
\item For every $\ep >0$, there are sums of squares of polynomials $s$ and $t$
in $\rx$ such that\label{ptentep}
\begin{equation}\label{ptentrep}
f+\ep=s+t(1-\|\nabla f\|^2\|\x\|^2).
\end{equation}
\end{enumerate}
\end{theorem}

\begin{proof}
First of all, the polynomial $g:=1-\|\nabla f\|^2\|\x\|^2$ is a polynomial
describing the principal gradient tentacle
$$S:=\{x\in\R^n\mid g(x)\ge 0\}=\ptent f.$$
Because sums of squares of polynomials are globally nonnegative on
$\R^n$, identity (\ref{ptentrep}) can be viewed as a certificate for
$f\ge -\ep$ on $S$. Hence it is clear that (\ref{tentep}) implies
(\ref{tentnn}). For the reverse implication, we apply Theorem \ref{schmuedgen}
(with $m=1$ and $g_1:=g$) to $f+\ep$ instead of $f$. We only have to check the
hypotheses. Condition (\ref{sbounded}) is clear from Lemma
\ref{lojasiewicz}. By Theorem \ref{parusinski}, we have that
$\ran(f,S)$ is a finite set if $f$ has only isolated singularities at infinity.
If $\ptent f$ is compact, the set $\ran(f,S)$ is even empty.
Since $f\ge 0$ on $S$ by hypothesis, this set
contains clearly only nonnegative numbers. This shows condition
(\ref{sfinite}), i.e.,
$\ran(f+\ep,S)=\ep+\ran(f,S)$ is a finite subset of $\R_{>0}$.
Finally, the hypothesis $f\ge 0$ on $S$ gives $f+\ep>0$ on $S$ which is
condition (\ref{spositive}). Therefore (\ref{ptentnn}) and (\ref{ptentep}) are
proved to be equivalent.
The equivalence of (\ref{pglobnn}) and (\ref{ptentnn}) is an immediate
consequence of Theorem \ref{critical}.
\end{proof}

\begin{remark}
Let $f\in\rx$ be bounded from below and $\ptent f$ be compact. Then $f$ attains
its infimum $f^*$. To see this, observe that the
equivalence of \eqref{pglobnn} and \eqref{ptentnn} in the preceding theorem
implies
\begin{align*}
f^*&=\sup\{a\in\R\mid f-a\ge 0\text{\ on $\R^n$}\}\\
&=\sup\{a\in\R\mid f-a\ge 0\text{\ on $\ptent f$}\}\\
&=\min\{f(x)\mid x\in\ptent f\}.
\end{align*}
\end{remark}

The following observation is proved in the same way than Remark \ref{homobs}.

\begin{remark}\label{homobs2}
If $f$ is a sum of squares in the ring $\R[[\x]]$ of formal power series, then
its lowest (non-vanishing) homogeneous part must be a sum of squares in
$\rx$.
\end{remark}

\begin{remark}\label{gen3}
There are polynomials $f\in\rx$ such that $f\ge 0$ on $\R^n$ but there is
no representation \eqref{ptentrep} for $\ep=0$. To see this, take a polynomial
$f\in\rx$ such that $f\ge 0$ on $\R^n$ but $f$ is not a sum of squares in the
ring $\R[[\x]]$ of formal power series (the Motzkin polynomial from
\eqref{motzkin} is such an example by the preceding remark). Then a
representation \eqref{ptentrep} with $\ep=0$ is impossible since the
polynomial $1-\|\nabla f\|^2\|\x\|^2$ has a positive constant term and
is therefore a square in $\R[[\x]]$.
\end{remark}

\subsection{Optimization using the gradient tentacle and sums of squares}
\label{psdp}

Theorem \ref{ptheorem} shows that under certain conditions, computation
of $f^*$ amounts to computing the supremum over all $a$ such that
$f-a=s+t(1-\|\nabla f\|^2\|\x\|^2)$ for some sums of squares $s$ and $t$
in $\rx$. As sketched in the introduction, sums of squares
\emph{of bounded degree} can be nicely parametrized by positive
semidefinite matrices. This motivates the following definition.

\begin{definition}\label{fkdef}
For all polynomials $f\in\rx$ and all $k\in\N_0$, we define
$f_k^*\in\R\cup\{\pm\infty\}$
as the supremum over all $a\in\R$ such that $f-a$ can be written as a sum
\begin{equation}\label{contraint}
f-a=s+t(1-\|\nabla f\|^2\|\x\|^2))
\end{equation}
where $s$ and $t$ are sums of squares of polynomials with $\deg t\le 2k$.
\end{definition}

Here and in the following, we use the convention that the degree of the zero
polynomial is $-\infty$ so that $t=0$ is allowed in the above definition.
Note that when the degree of $t$ in (\ref{contraint}) is restricted then
automatically also the degree of $s$.

Therefore the problem of computing $f_k^*$ can be written as an SDP. How to
do this, is already suggested in our introduction. It goes exactly like in the
well-known method of
Lasserre for optimization of polynomials on compact basic closed semialgebraic
sets. We refer to \cite{l1,m1,sr2} for the details. There are anyway several
toolboxes for MATLAB (a software for numerical computation) which can be used
to create and solve the corresponding SDPs without knowing these details.
The toolboxes we know are YALMIP \cite{löf} (which is very flexible and
good for much more than sums of squares stuff), SOSTOOLS \cite{pps} (which has
a very flexible and nice syntax), GloptiPoly \cite{hl} (very easy to use for
simple problems) and SparsePOP \cite{kkw} (specialized for sparse polynomials).
Besides MATLAB and such a toolbox one needs also an SDP solver for which the
toolbox provides an interface.

A side remark that we want to make here is that to each SDP there is a dual
SDP and it is desirable from the theoretical and practical point of view that
\emph{strong duality} holds, i.e., the optimal value of the primal
and dual SDP coincide. For the SDPs arising from Definition \ref{fkdef},
strong duality holds. This follows from the fact that principal gradient
tentacles (unlike gradient varieties) always have non-empty interior (they
always contain a small neighbourhood of the origin). For a proof confer
\cite[Theorem 4.2]{l1}, \cite[Corollary 3.2]{m1} or \cite[Corollary 21]{sr2}.
Here we will not define the dual SDP nor discuss its interpretation in terms
of the so-called moment problem.

Recalling the definition of $f^\sos$ in (\ref{sos}), we have obviously
\begin{equation}\label{chain}
f^\sos\le f_0^*\le f_1^*\le f_2^*\le\dots
\end{equation}
and if $f$ is bounded from below, then all $f_k^*$ are lower bounds
(perhaps $-\infty$) of $f^*$ by Theorem \ref{critical}. Note that the technique
from Jibetean and Laurent (see Subsection \ref{jlsubsection} above) gives
upper bounds for $f^*$ so that it complements nicely our method.
It is easy to see that Theorem \ref{ptheorem} can be expressed in terms of the
sequence $f_0^*,f_1^*,f_2^*,\dots$ as follows. 

\begin{theorem}\label{popt}
Let $f\in\rx$ be bounded from below. Suppose that $f$ has only isolated
singularities at infinity (e.g., $n=2$) or the principle gradient tentacle
$S(\nabla f)$ is compact. Then the sequence
$(f_k^*)_{k\in\N}$ converges monotonically increasing to $f^*$.
\end{theorem}

The following example shows that it is unfortunately in general not true that
$f_k^*=f^*$ for big $k\in\N$.

\begin{example}
Let $f$ be the Motzkin polynomial from \eqref{motzkin}. By Theorem \ref{popt},
we have
$\lim_{k\to\infty}f_k=0$. But it is not true that $f_k=0$ for some $k\in\N$.
By Definition \ref{fkdef}, this would imply that for all $\ep>0$, there is an
identity \eqref{ptentrep} with sums of squares $s$ and $t$ such that
$\deg s\le k$. Because $\ptent f$ has non-empty interior (note that
$\nabla f(1,1,1)=0$ since $f(1,1,1)=0$), we can use
\cite[Proposition 2.6(b)]{ps} (see \cite[Theorem 4.5]{sr2} for a
more elementary exposition) to see that such an identity would then also have
to exist for $\ep=0$. But this is impossible as we have seen in Remark
\ref{gen3}.
\end{example}

Unfortunately, the assumption that $f$ is bounded from below
is necessary in Theorem \ref{popt} as shown by the following trivial
example.

\begin{example}\label{trivex}
Consider $f:=X\in\R[X]$ (i.e., let $n=1$ and write $X$ instead of $X_1$).
Then $K(f)=\emptyset$, $\ptent f=[-1,1]$ and $(f_k^*)_{k\in\N}$ converges
monotonically increasing to
$\inf\{f(x)\mid -1\le x\le 1\}=-1\neq -\infty=f^*$.
\end{example}

\begin{openproblem}\label{open}
Do Theorems \ref{ptheorem} and \ref{popt} hold without the
hypothesis that $f$ has only isolated singularities at infinity or
$S(\nabla f)$ is compact?
\end{openproblem}

By the above arguments, it is easy to see that this question could be answered
in the affirmative if $\ran(f,\ptent f)$ were finite for all polynomials
$f\in\rx$ bounded from below on $\R^n$. But this is not true as the following
counterexample shows. We are grateful to Zbigniew Jelonek for pointing out to
us this adaption of an example of Parusi\'nski
\cite[Example 1.11]{p2}.

\begin{example}\label{neediso}
Consider the polynomial
$h:=X+X^2Y+X^4YZ\in\R[X,Y,Z]$,
set $f:=h^2$ and define for fixed $a>0$ the curve
$$\ga:\pos\to\R^3:s\mapsto\left(s,\frac{2a}{s^2},-\frac{\left(1+\frac s{4a}
  \right)}{2s^2}\right).$$
Observe that
$$h(\ga(s))=\frac 34s+a\qquad\text{and}\qquad
  \frac{\partial h}{\partial X}(\ga(s))=0$$
and therefore $f(\ga(s))=(\frac 34s+a)^2$ and
$$\|\nabla f\|^2(\ga(s))=4f\|\nabla h\|^2(\ga(s))=
  4s^4\left(\frac 34s+a\right)^2\left(\left(\frac 12-\frac s{8a}\right)^2+
  (2a)^2\right).$$
It follows that $\|\nabla f\|^2(\ga(s))\|\ga(s)\|^2$ equals
$$\left(4s^6+16a^2+\left(1+\frac s{4a}\right)^2\right)
  \left(\frac 34s+a\right)^2\left(\left(\frac 12-\frac s{8a}\right)^2+
  (2a)^2\right)$$
which tends to $(16a^2+1)a^2(1/4+4a^2)$ for $s\to 0$.
We now see that for $s\to 0$, $\|\ga(s)\|$ tends to infinity, $f(\ga(s))$ tends
to $a^2$ and, when $a$ is a sufficiently small positive number,
$\|\nabla f\|^2(\ga(s))\|\ga(s)\|^2$ tends to a real number smaller than $1$.
This shows that $a^2\in\ran(f,\ptent f)$ for every sufficiently small positive
number $a$. Hence $f$ is an example of a polynomial bounded from below such
that $\ran(f,\ptent f)$ is infinite.
\end{example}

\subsection{Implementation in YALMIP}\label{yalmip}

We show here how to encode computation of $f_k^*$ (as well as
of $f_{-1}^*:=f^\sos$) for any $k\in\N$ with YALMIP. First you have
to declare the variables appearing in the polynomial
$f$ (here $x$ and $y$) as well as the variable $a$ to maximize.
\begin{verbatim}
sdpvar x y a
\end{verbatim}
Now you specify the polynomial $f$ and the degree bound $k$ ($-1$ for
computing $f^\sos$). Here we take the dehomogenization $f:=M(X,Y,1)$ where
$M$ is the Motzkin polynomial introduced in (\ref{motzkin}).
\begin{verbatim}
f = x^4 * y^2 + x^2 * y^4 - 3 * x^2 * y^2 + 1, k = 0
\end{verbatim}
Now compute the partial derivatives with respect to the variables (here $x$ and
$y$) and specify the polynomial $g$ defining the gradient tentacle.
\begin{verbatim}
df = jacobian(f, [x y]), g = 1 - (df(1)^2 + df(2)^2) * (x^2 + y^2)
\end{verbatim}
Define a polynomial variable $t$ of degree $\le 2k$ and impose the constraints
that $t$ and $f-a-tg$ are sums of squares (for some reason the current
version of YALMIP does here not
accept a degree zero polynomial $t$ so that this has to be modeled as a scalar
variable).
\begin{verbatim}
if k > 0
  v = monolist([x; y], 2*k), coeffVec = sdpvar(length(v), 1)
  t = coeffVec' * v
  constraints = set(sos(f - a - t * g)) + set(sos(t))
elseif k == 0
  coeffVec = sdpvar(1, 1), t = coeffVec
  constraints = set(sos(f - a - t * g)) + set(t > 0)
else
  coeffVec = []
  constraints = set(sos(f - a))
end
\end{verbatim}
Now solve the SDP and output the result for $a$.
\begin{verbatim}
solvesos(constraints, -a, [], [a; coeffVec]), double(a)
\end{verbatim}

\subsection{Implementation in SOSTOOLS}\label{sostools}

Below we give an SOSTOOLS code which even slightly easier to read but
essentially analogous to the YALMIP code. In contrast to the YALMIP code
above, the MATLAB Symbolic Math Toolbox is required to execute the code below.
\begin{verbatim}
syms x y a t
f = x^4 * y^2 + x^2 * y^4 - 3 * x^2 * y^2 + 1, k = 0
df = jacobian(f, [x y]), g = 1 - (df(1)^2 + df(2)^2) * (x^2 + y^2)
prog = sosprogram([x; y], a)

if k > 0
  v = monomials([x; y], [0 : k]), [prog, t] = sossosvar(prog, v)
  prog = sosineq(prog, f - a - t * g)
elseif k == 0
  prog = sosdecvar(prog, t), prog = sosineq(prog, t)
  prog = sosineq(prog, f - a - t * g)
else
  prog = sosineq(prog, f - a)
end

prog = sossetobj(prog, -a), prog = sossolve(prog)
sosgetsol(prog, a)
\end{verbatim}

\subsection{Numerical results}\label{numerical}
The following examples have been computed on an ordinary PC with MATLAB 7,
YALMIP 3 and the SDP solver SeDuMi 1.1. Most of the computations took a
few seconds, some of them a few minutes. The first example corresponds
exactly to the code in Subsection \ref{yalmip}. To compute the others,
the variables, the polynomial $f$ and the degree bound $k$ has to be changed
in that code.

\begin{example}\label{exmotzkinxy}
Let $f:=M(X,Y,1)$ be the dehomogenization of the Motzkin polynomial $M$
from (\ref{motzkin}), i.e.,
$f:=M(X,Y,1)=X^4Y^2+X^2Y^4-3X^2Y^2+1\in\R[X,Y]$. We have $f^*=0$ but
$f^\text{sos}=-\infty$ (the latter is an easy exercise). If we execute the
program from Subsection \ref{yalmip}
with $k=-1$ instead of $k=0$, the computer answers that the SDP is infeasible
which means indeed that $f^\text{sos}=-\infty$. Executing the same program
for $k=0,1,2$ yields $f_0^*\approx -0.0017$, $f_1^*\approx-0.0013$ and
$f_2^*\approx 0.000066$ which is already very close to $f^*=0$. By Theorem
\ref{popt}, the sequence $f_0,f_1,f_2,\dots$ converges monotonically to
$f^*=0$. But the computed value $f_2^*\approx 0.000066$ is positive so that
there are obviously numerical problems. Confer \cite[Example 2]{ps}.
\end{example}

\begin{example}\label{exmotzkinxz}
Define $f:=M(X,1,Z)\in\R[X,Z]$ where $M$ is the Motzkin polynomial
from (\ref{motzkin}), i.e., $f=X^4+X^2+Z^6-3X^2Z^2\in\R[X,Z]$. Computation
yields $f^\sos\approx -0.1780$, $f_0^*\approx -5.1749\cdot 10^{-5}$,
$f_1^*\approx-1.2520\cdot 10^{-7}$ and $f_2^*=8.7662\cdot 10^{-10}$ which
``equals numerically'' $f^*=0$. This is in accordance with Theorem
\ref{ptheorem} which guarantees convergence to $f^*$ since we are in the two
variable case. Confer \cite[Example 3]{ps}.
\end{example}

\begin{example}\label{exberg}
Consider the Berg polynomial
$f:=X^2Y^2(X^2+Y^2-1)\in\R[X,Y]$
with global minimum $f^*=-1/27$ attained in $(\pm 1/\sqrt 3,\pm 1/\sqrt 3)$.
We have $f^\sos=-\infty$ and running the corresponding program gives indeed
an output saying that the corresponding SDP is infeasible. The computed optimal
values of the first principal tentacle relaxations are
$f_0^*\approx -0.0564$, $f_1^*\approx -0.0555$,
$f_2^*\approx -0.0371$ and $f_3^*\approx -0.0370\approx -1/27=f^*$.
Confer \cite[Example 3]{l1}, \cite[Example 3]{nds} and
\cite[Example 4]{jl}.
\end{example}

\begin{example}\label{ex11}
Being a polynomial in two variables of degree at most four,
we have that for $f:=(X^2 + 1)^2 + (Y^2 + 1)^2 - 2(X + Y + 1)^2\in\R[X,Y]$,
$f-f^*$ must be a sum of squares (see introduction) whence $f^*=f^\sos$.
By computation, we obtain for all values $f^\sos, f^*_0, f^*_1, f^*_2$
approximately $-11.4581$.  That all these computed values are the same can be
expected by $f^*=f^\sos$ and the monotonicity (\ref{chain}).
Confer \cite[Example 2]{l1} and \cite[Example 3]{jl}.
\end{example}

\begin{example}\label{exlax}
In \cite{ll}, it is shown that
\begin{align*}
f:=\sum_{i=1}^5\prod_{j\neq i}(X_i-X_j)\in\R[X_1,X_2,X_3,X_4,X_5]
\end{align*}
is nonnegative on $\R^5$ but not a sum of squares of polynomials. Therefore
$f^\sos=-\infty$ by Remark \ref{homobs} since $f$ is homogeneous.
The SDP solver detects indeed infeasibility of the corresponding SDP.
We have computed $f^*_0\approx-0.2367$, $f^*_1\approx -0.0999$ and
$f^*_2\approx -0.0224$. Solving the SDP relaxation computing $f^*_2$
took already the time of a coffee break. As in \cite[Example 6]{jl}, we observe
therefore that minimizing $f$ is after the change of variables
$X_{i}\mapsto X_1-Y_i$ ($i=2,3,4,5$) equivalent to minimizing
$$h:=Y_2Y_3Y_4Y_5+\sum_{i=2}^5(-Y_i)
     \prod_{j\neq i}(Y_j-Y_i)\in\R[Y_2,Y_3,Y_4,Y_5].$$
Computing $h^\sos$ results in infeasibility. The numerical results using the
principle gradient tentacle are $h_0^*\approx -0.2380$,
$h_1^*\approx -0.0351$, $h_2^*\approx -0.0072$, $h_3^*\approx -0.0019$ and
$h_4^*\approx -0.00086285$ which is already very close to $h^*=0$.
The condition in Theorem \ref{popt} is satisfied neither for $f$ nor
for $h$ and yet it seems that we have convergence to $h^*$. This is a typical
observation that might give hope that Open Problem \ref{open} has a positive
answer.
\end{example}

\begin{example}\label{exmain}
Consider once more the polynomial $f = (1-XY)^2+Y^2$ from (\ref{doesnot})
and Example \ref{bex}
that does not attain its infimum $f^*=0$ on $\R^2$. Since this polynomial is
by definition a sum of squares, we have $f^\sos=0=f^*$ and therefore
$f_k^*=0$ for all $k\in\N$ by (\ref{chain}). By computation, we get
$f^\sos\approx 1.5142\cdot 10^{-12}$ which is almost zero but also
$f_0^*\approx 0.0016$, $f_1^*\approx 0.0727$ and $f_2^*\approx 0.1317$
which shows that there are big numerical problems. We have verified that the
corresponding SDPs have nevertheless been solved quite accurately. The
problem is that small numerical errors in the coefficients of a polynomial
can perturb its infimum quite a lot whenever the infimum is not attained (or
attained very far from the origin). It should be subject to further research
how to fight this problem. Anyway, the gradient tentacle method still performs
in this example much better than the gradient
variety method which yields the wrong answer $1$ (as described in Subsection
\ref{gradvariety} above). The method of Jibetean and Laurent gives the
best results in this case \cite[Example 5]{jl}.
\end{example}

\subsection{Numerical stability}\label{stability}

If the coefficients of $f$ and $\|\nabla f\|\|\x\|$ have an order of magnitude
very different from $1$, then the defining polynomial
$g=1-\|\nabla f\|^2\|\x\|^2$
for the gradient tentacle should be better exchanged by
$R-\|\nabla f\|^2\|\x\|^2$
where $R$ is a real number of that order of magnitude. This is justified by
Remark \ref{nierem} above.

Example \ref{exmain} and other experiments that we did with polynomials bounded
from below that do not attain a minimum are a bit disappointing and show that
for this ``hard'' class of polynomials (exactly the class we were attacking), a
lot of work remains to be done, at least on the numerical side. The
corresponding semidefinite programs tend to be numerically unstable.

For polynomials attaining their minimum, the method in \cite{nds} is often
much more efficient, e.g., for Example \ref{exlax}.

\section{Higher gradient tentacles}\label{highersection}

In this section, we associate to every polynomial $f\in\rx$ a
\emph{sequence} of gradient tentacles. Each of these is defined by a
polynomial inequality just as the principal tentacle from Section
\ref{parusinskisection} was. But the degree of this polynomial inequality
for the $N$-th tentacle in this sequence will be roughly $2N$ times the
degree of $f$. This has the disadvantage that the corresponding SDP
relaxations get very big for large $N$. Also, we have to deal for
\emph{each} $N$ with a sequence of SDPs. All in all, we have therefore a
double sequence of SDPs. The advantage is however that we can prove a sums of
squares representation theorem (Theorem \ref{zeltschmid}) applicable for
\emph{all} $f\in\rx$ bounded from below independently of what is the answer
to Open Problem \ref{open}.
Again, we think that this theorem is also of theoretical
interest. Implementation of the higher gradient tentacle method is analogous
to Subsections \ref{yalmip} and \ref{sostools}. This time we do not give
numerical examples because of Open Problem \ref{open}, Remark \ref{generic}
and numerical problems for big $N$.

\begin{definition}\label{tentdef}
For $f\in\rx$ and $N\in\N$, we call
$$\tent fN:=\{x\in\R^n\mid\|\nabla f(x)\|^{2N}(1+\|x\|^2)^{N+1}\le 1\}$$
the $N$-th \emph{gradient tentacle} of $f$.
\end{definition}

A trivial fact that one should keep in mind is that
$\|\nabla f(x)\|^2(1+\|x\|^2)\le 1$ and in particular
$\|\nabla f(x)\|\|x\|\le 1$ for all $x\in\tent fN$. This shows that
$$V(\nabla f)\cap\R^n\subseteq\tent f1\subseteq\tent f2\subseteq\tent
  f3\subseteq\ldots\subseteq\ptent f.$$
The definition of $\tent fN$ is motivated by the following definition which
is taken from \cite[page 79]{kos}.

\begin{definition}\label{critvalinfndef}
Suppose $f\in\rx$ and $N\in\N$. The set $K_\infty^N(f)$ consists of all
$y\in\R$ for which there exists a sequence $(x_k)_{k\in\N}$ in $\R^n$ such
that
\begin{equation}\label{gradrn}
\lim_{k\to\infty}||x_k||=\infty,\qquad
\lim_{k\to\infty}\|\nabla f(x_k)\|\|x_k\|^{1+\frac 1N}=0\quad\text{and}\qquad
\lim_{k\to\infty}f(x_k)=y.
\end{equation}
\end{definition}

Clearly, we have
$$K_\infty^1(f)\subseteq K_\infty^2(f)\subseteq K_\infty^3(f)\subseteq\dots
  \subseteq K_\infty(f).$$
The next lemma says that this chain actually gets stationary and reaches
$K_\infty(f)$. For the proof, we refer to \cite[Lemma 3.1]{kos}.

\begin{lemma}[Kurdyka, Orro and Simon]\label{koslemma}
For all $f\in\rx$, there exists $N\in\N$ such that
$$K_\infty(f)=K_\infty^N(f).$$
\end{lemma}

Now we prove for sufficiently large gradient tentacles what was Corollary
\ref{critical} for the principal gradient tentacle (which contains all
higher gradient tentacles).

\begin{theorem}\label{mohab}
Suppose $f\in\rx$ is bounded from below.
Then $f^*\in K(f)$ and there is $N_0\in\N$ such that for all
$N\ge N_0$,
\begin{equation}\label{mohabeq}
f^*=\inf\{f(x)\mid x\in\tent fN\}.
\end{equation}
\end{theorem}

\begin{proof}
We know already from Theorem \ref{critical} that $f^*\in K(f)$. By
Proposition \ref{union}, at least one of the following two cases therefore must
occur. The first case is that $f^*\in K_0(f)$. Then $f^*$ is attained by $f$
on its gradient variety and therefore on the $N$-th gradient tentacle for
actually all $N\in\N$. Hence we can set $N_0:=1$. In the second case
$f^*\in K_\infty(f)$, we can choose some $N_0\in\N$ such that
$f^*\in K_\infty^N(f)$ by the previous Lemma. Then
$f^*\in K_\infty^N(f)$ for any $N\ge N_0$. This means that there exists
a sequence $(x_k)_{k\in\N}$ satisfying (\ref{gradrn}).
Therefore $\|\nabla f(x)\|\|x_k\|^{1+1/N}\le\frac 12$ and consequently
$$\|\nabla f(x_k)\|^{2N}(1+\|x_k\|^2)^{N+1}\le
  \|\nabla f(x_k)\|^{2N}(2\|x_k\|^2)^{N+1}\le 1
$$
for all large $k$ since $\|x_k\|\ge 1$ and $2^{N+1}\le 2^{2N}$.
This shows that $x_k\in\tent fN$ for all large $k$ which implies our
claim.
\end{proof}

The great advantage of the higher gradient tentacles over the principal one is
that they are \emph{always} small enough to admit only finitely many
asymptotic values, i.e., there is no counterpart to Example \ref{neediso}.

\begin{theorem}\label{finiteasymptotic}
For every $f\in\rx$, $\ran(f,\ptent f)\subseteq K_\infty(f)$. In particular,
every $f\in\rx$ has only finitely many asymptotic values on each
of its higher gradient tentacles, i.e., the set
$\ran(f,\tent fN)$ is finite for all $N\in\N$.
\end{theorem}

\begin{proof}
Let $y\in\R$ be such that (\ref{asymptotic}) holds for some sequence
$(x_k)_{k\in\N}$ of points $x_k\in\tent fN$. By Definition \ref{tentdef},
$$\|\nabla f(x_k)\|^N\|x_k\|^N\le\frac 1{\|x_k\|}\to 0\qquad
  \text{for $k\to\infty$}$$
implying (\ref{gradr}). This shows $y\in K_\infty(f)$.
\end{proof}

\subsection{Higher gradient tentacles and sums of squares}

We are now able to prove the third important sums of squares representation
theorem of this article besides Theorems \ref{schmuedgen} and \ref{ptheorem}.

\begin{theorem}\label{zeltschmid}
For all $f\in\rx$ bounded from below, there is $N_0\in\N$ such that for all
$N\ge N_0$, the following are equivalent.
\begin{enumerate}[(i)]
\item $f\ge 0$ on $\R^n$\label{globnn}
\item $f\ge 0$ on $\tent fN$\label{tentnn}
\item For every $\ep >0$, there are sums of squares of polynomials $s$ and $t$
in $\rx$ such that\label{tentep}
\begin{equation}\label{tentrep}
f+\ep=s+t(1-\|\nabla f\|^{2N}(1+\|\x\|^2)^{N+1}).
\end{equation}
\end{enumerate}
Moreover, these conditions are equivalent for all $f$ attaining a minimum on
$\R^n$ and all $N\in\N$. Finally, (\ref{tentnn}) and (\ref{tentep}) are
equivalent for all  $f\in\rx$ and $N\in\N$.
\end{theorem}

\begin{proof}
We first show that (\ref{tentnn}) and (\ref{tentep}) are always equivalent.
To see this, observe that $g_1:=1-\|\nabla f\|^{2N}\|\x\|^{2N+2}$ is a
polynomial that defines the set $S:=\{x\in\R^n\mid g_1\ge 0\}=\tent fN$.
Because sums of squares of polynomials are globally nonnegative on
$\R^n$, identity (\ref{tentrep}) can be viewed as a certificate for
$f\ge -\ep$ on $S$. Hence it is clear that (\ref{tentep}) implies
(\ref{tentnn}). For the reverse implication, we apply Theorem
\ref{schmuedgen} to $f+\ep$ instead of $f$. We only have to check the
hypotheses. Condition (\ref{sbounded}) is clear from Lemma
\ref{lojasiewicz}. By Corollary \ref{finiteasymptotic}, we have that
$\ran(f,S)$ is a finite set. Since $f\ge 0$ on $S$ by hypothesis, this set
contains clearly only nonnegative numbers. This shows condition
(\ref{sfinite}), i.e.,
$\ran(f+\ep,S)=\ep+\ran(f,S)$ is a finite subset of $\R_{>0}$.
Finally, the hypothesis $f\ge 0$ on $S$ gives $f+\ep>0$ on $S$ which is
condition (\ref{spositive}).

Now suppose that $f\in\rx$ attains a minimum $f(x^*)=f^*$ in a point
$x^*\in\R^n$. Then $\nabla f(x^*)=0$ and therefore $x^*\in\tent fN$ for all
$N\in\N$. This shows that (\ref{globnn}) and (\ref{tentnn}) are in this
case equivalent for \emph{all} $N\in\N$.

 By what has already been proved,
it remains only to show that (\ref{globnn}) and (\ref{tentnn}) are equivalent
for large $N\in\N$ when $f\in\rx$ is bounded from below but does not attain a
minimum. But in this
case, (\ref{mohabeq}) holds by Theorem \ref{mohab} yielding the equivalence
of the first two conditions.
\end{proof}

Without needing it for our application, we draw the following immediate
corollary.
Taking $N=1$ in the second part of this corollary yields Theorem \ref{nie}
above of Nie, Demmel and Sturmfels.

\begin{corollary}\label{gradpos}
Suppose $f\in\rx$ and $f\ge 0$ on $V(\nabla f)\cap\R^n$. Then $f+\ep$ is for
all $\ep>0$ a sum of squares modulo any principal ideal generated by a power
of the polynomial
$\|\nabla f\|^2(1+\|\x\|^2)$, i.e., for every $\ep>0$ and $N\in\N$,
there is a sum of squares $s$ in $\rx$ and a polynomial $p\in\rx$ such that
$$f=s+p(\|\nabla f\|^2(1+\|\x\|^2))^N.$$
In particular, $f+\ep$ is for all $\ep>0$ a sum of squares modulo each power
of its gradient ideal, i.e., for every $\ep>0$ and $N\in\N$, there is a sum
of squares $s$ in $\rx$ such that $$f\in s+(\nabla f)^N.$$
\end{corollary}

\begin{proof}
The second claim follows from the first one.
The first claim follows immediately from
implication (\ref{globnn})$\implies$(\ref{tentep}) in Theorem
\ref{zeltschmid} which always holds for all $N\in\N$.
\end{proof}

\subsection{Optimization using higher gradient tentacles and sums of squares}

The following definition can be motivated in the same way than Definition
\ref{fkdef} in Section \ref{parusinskisection}.

\begin{definition}\label{fnkdef}
For all polynomials $f\in\rx$, all $N\in\N$ and all $k\in\N_0$, we define
$f_{N,k}^*\in\R\cup\{\pm\infty\}$
as the supremum over all $a\in\R$ such that $f-a$ can be written as a sum
\begin{equation}\label{constr}
f-a=s+t(1-\|\nabla f\|^{2N}(1+\|\x\|^2)^{N+1})
\end{equation}
where $s$ and $t$ are sums of squares of polynomials with $\deg t\le 2k$.
\end{definition}

Again, like in Section \ref{parusinskisection} outlined, computation of
$f_{N,k}$ amounts to solving an SDP for each fixed $N\in\N$ and $k\in\N_0$.
Recalling the definition of $f^\sos$ in (\ref{sos}), we have for each fixed
$N\in\N$,
$$f^\sos\le f_{N,0}^*\le f_{N,1}^*\le f_{N,2}^*\le\dots$$
and if $f$ is bounded from below, then all $f_{N,k}^*$ are lower bounds
of $f^*$ by Theorem \ref{mohab}.
It would be desirable to have also information
how the $f_{N,k}$ are related to each other when not only $k$ but also
$N$ varies. All we know about that is the following proposition.

\begin{proposition}\label{tentinccert}
For all $f\in\rx$, $N\in\N$ and $k\in\N_0$,
$$f^*_{N+1,k}\le f^*_{N,k+d}.$$
\end{proposition}

\begin{proof}
Let us define the polynomials $h_N$ like in (\ref{hndef}) and substitute in the
identity proved in Lemma \ref{hnlemma}, the polynomials $\|\nabla f\|^2$ for
$Y$ and $\|\x\|^2$ for $\x$. Then we get
\begin{equation}\label{subst}
1-\|\nabla f\|^{2(N+1)}(1+\|\x\|^2)^{N+2}
=p+q(1-\|\nabla f\|^{2N}(1+\|\x\|^2)^{N+1}).
\end{equation}
where $p$ and
$$q:=\left(1+\frac 1N\right)\|\nabla f\|^2(1+\|\x\|^2)$$
are sums of squares of polynomials. The degree of $q$ is no higher than
$2(d-1)+2=2d$. Now if for $a\in\R$ we have an identity
$$f-a=s+t(1-\|\nabla f\|^{2(N+1)}(1+\|\x\|^2)^{N+2})$$
for sums of squares $s$ and $t$ with $\deg t\le 2k$, then for the same $a$
$$f-a=(s+tp)+tq(1-\|\nabla f\|^{2N}(1+\|\x\|^2)^{N+1})$$
and $\deg(tq)\le 2(k+d)$.
\end{proof}

We conclude by interpreting Theorem \ref{zeltschmid} as a convergence result
concerning the optimal values $f_{N,k}^*$ of the proposed relaxations.
This is the counterpart to Theorem \ref{popt} from Section \ref{sossection}.

\begin{theorem}\label{opt}
For all $f\in\rx$ bounded from below, $(f_{N,k}^*)_{k\in\N}$ converges
monotonically increasing to $f^*$ provided that $N\in\N$ is sufficiently large
(depending on $f$). If $f$ attains a minimum on $\R^n$, $(f_{N,k}^*)_{k\in\N}$
converges monotonically increasing to $f^*$ no matter what $N\in\N$ is.
\end{theorem}

\section{Conclusions}\label{conclusion}
We have proposed a method for computing numerically the infimum of a real
polynomial in $n$ variables which is bounded from below on $\R^n$. Like in
\cite{jl} and \cite{nds}, the approach is to find semidefinite relaxations
relying on sums of squares certificates and critical point theory. As one could
expect, polynomials that do not attain a minimum on $\R^n$ (that are either
unbounded from below or have a finite infimum that is not attained) are
particularly hard to handle.
In \cite{jl}, this problem (among others) was solved by perturbing
the coefficients of the polynomial to guarantee a minimum (in particular,
boundedness from below). Though the results in \cite{jl} are quite good,
we are convinced that one should also look for other methods that avoid
perturbations and the danger of numerical ill-conditioning coming along with
them. Proving sums of squares representations for polynomials positive on their
gradient variety, it was shown by Nie, Demmel and Sturmfels \cite{nds} that an
approach without perturbation is possible. The computational performance of
their method is extremely good. However, for polynomials that do not attain a
minimum, their method yields wrong answers. Combining
considerable machinery from differential geometry and real algebraic geometry,
we have shown that part of this limitation can be removed. By using our
gradient tentacles instead of the gradient variety, polynomials that do not
attain a minimum but are bounded from below can also be handled. Our method
has three major problems. First, we do not address the important question of
how to check efficiently if a polynomial is bounded from below. For such
polynomials, our method still gives a wrong answer (see Example \ref{trivex}).
Second, it turns out that
solving semidefinite programs that arise from a polynomial that does not attain
a minimum takes sometimes surprisingly long time. And third, small numerical
inaccuracies might lead to big changes in the infimum of a polynomial if the
infimum is not attained. All three problems should be subject to further
research. Polynomials not attaining a minimum remain hard to handle in
practice. On the theoretical side, we have combined the theory of generalized
critical values with the the theory of real holomorphy rings and have obtained
new interesting characterizations of nonnegative polynomials.

\section*{Acknowledgments}
We are most grateful to Zbigniew Jelonek for the discussions in Passau where
he showed us Example \ref{neediso} and Parusinski's Theorem \ref{parusinski}.
Our thanks go also to Mohab Safey El Din for shifting our
attention to Theorem \ref{critical}, to Richard Leroy for helping to prove
Lemma \ref{hnlemma} and 
to Krzysztof Kurdyka for interesting discussions in Paris.


\begin{thebibliography}{KKW}

\bibitem[Ble]{ble} G. Blekherman.
There are Significantly More Nonnegative Polynomials than Sums of Squares.
Preprint.\\
\url{http://arxiv.org/abs/math.AG/0309130}

\bibitem[Cas]{cas} G. Cassier.
Probl\`eme des moments sur un compact de $\R^n$ et d\'ecomposition de
polyn\^omes \`a plusieurs variables.
J. Funct. Anal. {\bf 58}, 254--266 (1984)

\bibitem[HL]{hl} D. Henrion and J. Lasserre.
GloptiPoly: Global optimization over polynomials with Matlab and SeDuMi.
ACM Trans. Math. Softw. {\bf 29}, No. 2, 165--194 (2003)

\bibitem[JL]{jl} D. Jibetean and M. Laurent.
Semidefinite approximations for global unconstrained polynomial optimization.
SIAM Journal on Optimization {\bf 16}, No. 2, 490--514 (2005)

\bibitem[KKW]{kkw} M. Kojima, S. Kim and H. Waki.
Sparsity in sums of squares of polynomials.\\
Math. Program. {\bf 103}, No. 1 (A), 45--62 (2005)

\bibitem[KOS]{kos} K. Kurdyka, P. Orro and S. Simon:
Semialgebraic Sard theorem for generalized critical values.
J. Differ. Geom. {\bf 56}, No.1, 67-92 (2000)

\bibitem[Kri]{kri} J. Krivine.
Anneaux preordonnes.
J. Anal. Math. {\bf 12}, 307--326 (1964)

\bibitem[Kus]{kus} A. Kushnirenko.
Poly\`edres de Newton et nombres de Milnor.
Invent. Math. {\bf 32}, 1--31 (1976)

\bibitem[L1]{l1} J. Lasserre.
Global optimization with polynomials and the problem of moments.
SIAM J. Optim. {\bf 11}, No. 3, 796--817 (2001)

\bibitem[L2]{l2} J. Lasserre.
A sum of squares approximation of nonnegative polynomials.
SIAM J. Optim. {\bf 16}, No. 3, 751--765 (2006) 

\bibitem[LL]{ll} A. Lax and P. Lax.
On sums of squares.
Linear Algebra Appl. {\bf 20}, 71--75 (1978)

\bibitem[LN]{ln} J. Lasserre and T. Netzer.
SOS approximations of nonnegative polynomials via simple high degree
perturbations. To appear in Math. Z.\\
\url{http://arxiv.org/abs/math.math.AG/0510456}

\bibitem[L\"of]{löf} J. L\"ofberg.
YALMIP: A MATLAB toolbox for rapid prototyping of optimization problems\\
\url{http://control.ee.ethz.ch/~joloef/yalmip.php}

\bibitem[M1]{m1} M. Marshall.
Optimization of polynomial functions.
Can. Math. Bull. {\bf 46}, No. 4, 575--587 (2003)

\bibitem[M2]{m2} M. Marshall.
Representations of non-negative polynomials, degree bounds and applications
to optimization.
Preprint.\\
\url{http://math.usask.ca/~marshall/}

\bibitem[NDS]{nds} J. Nie, J. Demmel, and B. Sturmfels.
Minimizing polynomials via sum of squares over the gradient ideal.
Math. Prog., Ser. A, {\bf 106}, No. 3, 587--606 (2006)

\bibitem[Net]{net} T. Netzer.
High degree perturbations of nonnegative polynomials.
Diplomarbeit, Universit\"at Konstanz (2005)\\
\url{http://www.math.uni-konstanz.de/~netzer/}

\bibitem[P1]{p1} A. Parusi\'nski.
On the bifurcation set of a complex polynomial with isolated singularities at
infinity.
Compos. Math. {\bf 97}, No. 3, 369--384 (1995)

\bibitem[P2]{p2} A. Parusi\'nski.
A note on singularities at infinity of complex polynomials.
Proceedings of the Banach Center symposium on differential geometry and
mathematical physics in Spring 1995.
Banach Cent. Publ. {\bf 39}, 131--141 (1997)

\bibitem[PD]{pd} A. Prestel, C. Delzell.
Positive polynomials.
Monographs in Mathematics. Springer, Berlin (2001)

\bibitem[PPS]{pps} S. Prajna, A. Papachristodoulou, P. Seiler and P. Parrilo.
SOSTOOLS and its control applications.
Lecture Notes in Control and Information Sciences {\bf 312}, 273--292 (2005)

\bibitem[PS]{ps} P. Parrilo and B. Sturmfels.
Minimizing polynomial functions.
Ser. Discrete Math. Theor. Comput. Sci. {\bf 60}, 83-99 (2003)

\bibitem[Rab]{rab} P. Rabier.
Ehresmann fibrations and Palais-Smale conditions for morphisms of
Finsler manifolds.
Ann. Math. (2) {\bf 146}, No. 3, 647--691 (1997)

\bibitem[Rez]{rez} B. Reznick.
Some concrete aspects of Hilbert's 17th problem.
Contemp. Math. {\bf 253}, 251--272 (2000)

\bibitem[Sch]{sch} K. Schm\"udgen.
The $K$-moment problem for compact semi-algebraic sets.
Math. Ann. {\bf 289}, No. 2, 203--206 (1991)

\bibitem[Scd]{scd} J. Schmid.
On the degree complexity of Hilbert's 17th problem and the Real
Nullstellensatz.
Habilitationsschrift, Universit\"at Dortmund (1998)

\bibitem[Sho]{sho} N. Shor.
Class of global minimum bounds of polynomial functions.
Cybernetics {\bf 23}, No. 6, 731--734 (1987)

\bibitem[Shu]{shu} E. Shustin.
Critical points of real polynomials, subdivisions of Newton polyhedra and
topology of real algebraic hypersurfaces.
Transl., Ser. 2, Am. Math. Soc. {\bf 173}, 203--223 (1996)

\bibitem[Spo]{spo} S. Spodzieja.
Lojasiewicz inequalities at infinity for the gradient of a polynomial.
Bull. Pol. Acad. Sci., Math. {\bf 50}, No. 3, 273--281 (2002)

\bibitem[Sr1]{sr1} M. Schweighofer.
Iterated rings of bounded elements and generalizations of Schm\"udgen's
Positivstellensatz.
J. Reine Angew. Math. {\bf 554}, 19--45 (2003).
Erratum available at\\
\url{http://arxiv.org/abs/math.AC/0510675}

\bibitem[Sr2]{sr2} M. Schweighofer.
Optimization of polynomials on compact semialgebraic sets.
SIAM Journal on Optimization {\bf 15}, No. 3, 805--825 (2005)

\bibitem[SS]{ss} N. Shor and P. Stetsyuk.
Modified $r$-algorithm to find the global minimum of polynomial functions.
Cybern. Syst. Anal. {\bf 33}, No. 4, 482--497 (1997)

\bibitem[Ste]{ste} G. Stengle.
A Nullstellensatz and a Positivstellensatz in semialgebraic geometry.
Math. Ann. {\bf 207}, 87--97 (1974)

\bibitem[Tho]{tho} R. Thom.
Ensembles et morphismes stratifies.
Bull. Am. Math. Soc. {\bf 75}, 240--284 (1969)

\bibitem[Tod]{tod} M. Todd. Semidefinite Optimization.
Acta Numerica {\bf 10}, 515-560 (2001)
\end{thebibliography}
\end{document}